\documentclass{amsart}

\usepackage{fancyhdr}
\usepackage{lastpage}

\newcommand{\bdis}{\begin{displaymath}}
\newcommand{\edis}{\end{displaymath}}
\newcommand{\be}{\begin{equation}}
\newcommand{\ee}{\end{equation}}

\newcommand{\mcal}{\mathcal}

\newtheorem{theorem}{Theorem}
\newtheorem{lemma}[]{Lemma}

\theoremstyle{definition}

\theoremstyle{remark}
\newtheorem{remark}[theorem]{Remark}

\numberwithin{equation}{section}

\begin{document}

\title{Jacob's ladders and the almost exact asymptotic representation of the Hardy-Littlewood integral}

\author{Jan Moser}

\address{Department of Mathematical Analysis and Numerical Mathematics, Comenius University, Mlynska Dolina M105, 842 48 Bratislava, SLOVAKIA}

\email{jan.mozer@fmph.uniba.sk}


\begin{abstract}
In this paper we introduce a nonlinear integral equation such that
the system of global solution to this equation represents a class
of a very narrow beam at $T\to\infty$ (an analogue to the laser
beam) and this sheaf of solutions leads to an almost-exact
representation of the Hardy-Littlewood integral. The accuracy of
our result is essentially better than the accuracy of related results of Balasubramanian, Heath-Brown and Ivic.
\end{abstract}

\maketitle

\section{Introduction} 
\label{intro}

Let us remind that Hardy and Littlewood started to study the following integral in 1918:

\be \label{1}
\int_0^T\left|\zeta\left( \frac{1}{2}+it\right)\right|^2{\rm d}t=\int_0^T Z^2(t){\rm d}t ,
\ee
where
\bdis
Z(t)=e^{i\vartheta(t)}\zeta\left( \frac{1}{2}+it\right),\
\vartheta(t)=-\frac{1}{2}t\ln(\pi)+\mbox{Im} \ln\left[ \Gamma\left(\frac{1}{4}+i\frac{t}{2}\right)\right],
\edis
and they have derived the following formula (see \cite{4}, page 122, 151-156)
\be \label{2}
\int_0^T Z^2(t){\rm d}t\sim T\ln(T),\ T\to\infty .
\ee
In this paper we show that except the asymptotic formula (\ref{2}) that posses an unbounded error there is an infinite family of other asymptotic
representations of the Hardy-Littlewood integral (\ref{1}). Each member of this family is an almost-exact representation of the given
integral. The proof of this will be based on properties of a kind of functions having some canonical properties on the set of zeroes of the
function $\zeta(1/2+it)$. \\

(A)\ Let us remind further that in 1928 Titchmarsh has discovered a new treatment to the integral (\ref{1}) by which the Titchmarsh-Kober-Atkinson
(TKA) formula:
\be \label{3}
\int_0^\infty Z^2(t)e^{-2\delta t}{\rm d}t=\frac{c-\ln(4\pi\delta)}{2\sin(\delta)}+\sum_{n=0}^N c_n\delta^n+\mcal{O}\left(\delta^{N+1}\right) ,
\ee
where $\delta\to 0$, $C$ is the Euler constant, $c_n$ are constant depending upon $N$,was derived. \\

The TKA formula has been published for the first time in 1951 in the fundamental monograph by Titchmarsh (see \cite{8}, \cite{1},\cite{6},\cite{7}).
It was thought for about 56 years that the TKA formula is a kind of curiosity (see \cite{5}, page 139).
However, in this paper we show that the TKA formula itself contains new kinds of principles. \\

(B)\ Namely, in this work we introduce a new class of curves, which are the solutions to the following nonlinear integral equation:
\be \label{4}
\int_0^{\mu[x(T)]}Z^2(t)e^{-\frac{2}{x(T)}t}{\rm d}t=\int_0^T Z^2(t){\rm d}t ,
\ee
where the class of functions $\{\mu\}$ is specified as:
$\mu\in C^\infty ([y_0,\infty))$ is a monotonically increasing (to $+\infty$) function and it
obeys $\mu(y)\geq 7y\ln(y)$ .
The following holds true: for any $\mu\in\{ \mu\}$ it exists just one solution to the equation (\ref{4}):
\bdis
\varphi(T)=\varphi_\mu(T),\ T\in [T_0,\infty),\ T_0=T_0[\varphi],\ \varphi(T)\to\infty \ \mbox{as}\ T\to\infty .
\edis
Let us denote by the symbol $\{\varphi\}$ the system of these solutions. The function $\varphi(T)$ is related to the zeroes of the Riemann zeta-function
on the critical line by the following way. Let $t=\gamma$ be a zero of the function $\zeta(1/2+it)$ of the order $n(\gamma)$, where
$n(\gamma)=\mcal{O}(\ln(\gamma))$, (see \cite{3}, page 178). Then the points $[\gamma,\varphi(\gamma)],\ \gamma>T_0$ (and only these points) are
the inflection points with the horizontal tangent. In more details, it holds true the following system of equations:
\be \label{5}
\varphi^{\prime}(\gamma)=\varphi^{\prime\prime}(\gamma)=\dots=\varphi^{(2n)}(\gamma)=0,\ \varphi^{(2n+1)}(\gamma)\not=0 ,
\ee
where $n=n(\gamma)$. \\

With respect to this property an element $\varphi\in\{ \varphi\}$ is to be named as the Jacob's ladder leading to $[+\infty,+\infty]$ (the rungs of the Jacob's ladder are the
segments of the curve $\varphi$ lying in the neighborhoods of the points $[\gamma,\varphi(\gamma)],\ \gamma>T_0[\varphi]$). Finally,
also the composition of the functions $G[\varphi(T)]$ is to be named the Jacob's ladder if the following conditions are fulfilled: $G\in C^\infty([y_0,+\infty))$,
$G$ grows to $+\infty$ and $G$ has a positive derivative everywhere. \\

Let us mention that the mapping (the operator)
\bdis
\hat{H}:\ \{\mu\}\to\{\varphi\}
\edis
can be named the $Z^2$-mapping of the functions of the class $\{\mu\}$.  \\

(C)\ Jacob's ladder implies the following results:

\begin{itemize}
\item[(a)] An almost exact asymptotic formula for the Hardy-Littlewood integral. Let us mention that our new formula makes more exact also the
leading term in (\ref{2}):
\bdis
\int_0^T Z^2(t){\rm d}t\sim \frac{\varphi(T)}{2}\ln\left(\frac{\varphi(T)}{2}\right),\ T\to +\infty,\ \forall \varphi\in\{\varphi\} ,
\edis
i.e. the leading term in this formula is also the Jacob's ladder, and therefore we can say that the leading term has a very fine structure.
\item[(b)] The system $\{\varphi\}$ has the property of an "infinitely close approach" of any two Jacob's ladders at $T\to +\infty$.
\item[(c)] Our new formula for the Hardy-Littlewood integral is stable with respect to the choice of the elements from some
subset $\{\varphi\}^\star$ in $\{\varphi\}$.
\end{itemize}

\section{Results}

The following theorem holds true:

\begin{theorem} {\rm
The TKA formula implies: \\
(A) \
\be \label{6}
\int_0^T Z^2(t){\rm d}t=F[\varphi(T)]+r[\varphi(T)],\ T\geq T_0[\varphi],
\ee
where
\be \label{7}
F[y]=\frac{y}{2}\ln\left(\frac{y}{2}\right)+(c-\ln(2\pi))\frac{y}{2}+c_0,\
r[\varphi(T)]=\mcal{O}\left\{\frac{\ln(\varphi(T))}{\varphi(T)}\right\}=\mcal{O}\left(\frac{\ln(T)}{T}\right) .
\ee
(B)\ For all $\varphi_1(T),\varphi_2(T)\in\{\varphi\}$ we have
\be \label{8}
\varphi_1(T)-\varphi_2(T)=\mcal{O}\left(\frac{1}{T}\right),\ T\geq\max\{ T_0[\varphi_1],T_0[\varphi_2]\}.
\ee
(C)\ If the set $\{\mu(y_0)\}$ is bounded then for all $\varphi(T)$ in $\{\varphi\}$ and for any fixed $\varphi_0(T)\in\{\varphi\}$:
\be \label{9}
\varphi_0(T)-\frac{A}{T}<\varphi(T)<\varphi_0(T)+\frac{A}{T},\ T\geq T_0=\sup_{\mu}\{ \mu(y_0)\}.
\ee
Let us mention that the constants in the $\mcal{O}$-symbols do not depend upon the choice of $\varphi(T)$.
}
\end{theorem}

Let us remind the Balasubramanian's formula (see \cite{2}):
\be \label{10}
\int_0^T Z^2(t){\rm d}t=T\ln(T)+(2c-1-\ln(2\pi))T+\mcal{O}\left( T^{1/3+\epsilon}\right),
\ee
and $\Omega$-theorem of Good (see \cite{3}):
\be \label{11}
\int_0^T Z^2(t){\rm d}t-T\ln(T)-(2c-1-\ln(2\pi))T=\Omega\left( T^{1/4}\right) .
\ee

\begin{remark} {\rm
Combining the formulae (\ref{6}), (\ref{7}) and (\ref{10}) one obtains that:
\begin{itemize}
\item Formula (\ref{10}) possesses quite large uncertainty since the deviation from the value of (\ref{1}) is given by
$R(T)=\mcal{O}(T^{1/3+\epsilon})$, and (see (\ref{11})) since
\bdis
\overline{\lim_{T\to\infty}}|R(T)|=+\infty ,
\edis
this cannot be removed.
\item Following (\ref{7}) we have
\bdis
\lim_{T\to\infty}r[\varphi(T)]=0 ,
\edis
and this means that formula (\ref{6}) seems to be almost exact.
\end{itemize}
}
\end{remark}

\begin{remark} {\rm
We have found a new fact that the leading term in the Hardy-Littlewood integral is a ladder (see (\ref{6})), i.e. it has a fine structure. There is
no analogue of this in the formulae (\ref{2}) or (\ref{10}).
}
\end{remark}

\begin{remark} {\rm
We say explicitly that:
\begin{itemize}
\item Formula (\ref{8}) contains a new effect, namely, any two Jacob's ladders approach each other at $T\to\infty$.
\item (\ref{9}) implies that the representations (\ref{6}) and (\ref{7}) of the integral of Hardy-Littlewood (\ref{1}) are stable under
the choice of $\varphi(T)\in\{\varphi\}$ in the case of a bounded set $\{\mu(y_0)\}$.
\end{itemize}
}
\end{remark}

\begin{remark} {\rm
Following the second part of Remark 3 the representations (\ref{6}) and (\ref{7}) of the Hardy-Littlewood integral (\ref{1}) is
\emph{microscopically unique} in sense that any two Jacob's ladders $\varphi_1,\varphi_2,\ T\geq T_0$ (see (\ref{9})) cannot be distinguished at $T\to\infty$.
}
\end{remark}

In the fifth part of this work we establish the relation between the Jacob's ladder and the prime-counting function $\pi(T)$:
\bdis
\pi(T)\sim\frac{1}{1-c}\left\{T-\frac{\varphi(T)}{2}\right\},  \ T\to\infty .
\edis

\section{Existence of the Jacob's ladder}

The following lemma holds true:

\begin{lemma} {\rm
Let $\mu(y)\in\{\mu\}$ be fixed. Then there exists an unique solution $\varphi(T),\ T\geq T_0[\varphi]$ to the integral equation (\ref{4}) that \
obeys the property (\ref{5}).
}
\end{lemma}

\begin{proof}
By the Bonnet's mean-value theorem in the case $y>0$ and $\mu>0$ we have
\be \label{12}
\int_0^\mu Z^2(t)e^{-\frac{2}{y}t}{\rm d}t=\int_0^MZ^2(t){\rm d}t,\ M>0 ,
\ee
where $e^{-\frac{2}{y}t},\ t\in[0,\mu]$ is decreasing and equals $1$ at $t=0$. \\

First of all we will show that the formula (\ref{12}) maps to any fixed $\mu>0$ just one $M>0$, since the case $M_1\not=M_2$ is impossible
because it would mean that
\be \label{13}
\int_{M_1}^{M_2}Z^2(t){\rm d}t=0 .
\ee
Let $\mu(y)\in\{\mu\}$. Then the following formula
\be \label{14}
\int_0^{\mu(y)}Z^2(t)e^{-\frac{2}{y}t}{\rm d}t=\int_0^{M(y)}Z^2(t){\rm d}t
\ee
defines a function $M(y)=M_\mu(y),\ y\geq y_0$. Let the symbol $\{ M\}$ denote the class of the images of the elements $\mu(y)\in\{\mu\}$.  \\

The function $M(y)$ is positive and increases to $+\infty$. In fact, let
\be \label{15}
\Phi(y)=\int_0^{\mu(y)}Z^2(t)e^{-\frac{2}{y}t}{\rm d}t .
\ee
Then
\be \label{16}
\Phi'(y)=\frac{2}{y^2}\int_0^{\mu(y)}tZ^2(t)e^{-\frac{2}{y}t}{\rm d}t+Z^2[\mu(y)]e^{-\frac{2}{y}\mu(y)}\frac{{\rm d}\mu(y)}{{\rm d}y}>0
\ee
(the first term is obviously positive and the second one is non-negative). Thus, we see that the function $\Phi(y),\ y\in[y_0,+\infty)$ is
increasing. Subsequently, for any $y,\Delta y>0,y\geq y_0$ one has
\bdis
0<\Phi(y+\Delta y)-\Phi(y)=\int_{M(y)}^{M(y+\Delta y)}Z^2(t){\rm d}t \quad \Rightarrow \quad
M(y+\Delta y)>M(y) .
\edis

The function $M(y),\ y\geq y_0$ is continuous. In fact, for any fixed $\hat{y}\in[y_0,+\infty)$ we have
\bdis
M_1(\hat{y})=\limsup_{y\to\hat{y}^+}M(y),\quad M_2(\hat{y})=\liminf_{y\to\hat{y}^+}M(y) .
\edis
And with respect to the continuity of the left-hand side of eq. (\ref{14}) we obtain (see (\ref{13})):
\bdis
\int_{M_1(\hat{y})}^{M_2(\hat{y})}Z^2(t){\rm d}t=0 \quad M_1(\hat{y})=M_2(\hat{y})=M(\hat{y}+0) .
\edis

The existence of $M(\hat{y}-0)$ and the identities: $M(\hat{y})=M(\hat{y}-0)=M(\hat{y}+0)$ can be shown by analogy.

The function $M(y),\ y\geq y_0$ obeys the following properties:
\begin{itemize}
\item[(a)] It has a continuous derivative at any point $y$ such that $M(y)\not=\gamma$, where $\gamma$ is a zero of the function
$\zeta(1/2+it)$.
\item[(b)] It has a derivative equal to $+\infty$ at any point $y$ such that $M(y)=\gamma$, $\gamma$ mentioned above.
\end{itemize}
These properties of the $M$ function can be proved as follows. By (\ref{14}) and (\ref{15}) we have
\be \label{17}
\frac{\Phi(y+\Delta y)-\Phi(y)}{\Delta y}=\frac{M(y+\Delta y)-M(y)}{\Delta y}Z^2
\left\{ M(y)+\theta\cdot\left[ M(y+\Delta y)-M(y)\right]\right\},
\ee
where $\theta\in(0,1)$. Using the fact that $\Phi'(y)>0$ (see (\ref{16})) we can deduce from (\ref{17}) that:
\begin{itemize}
\item[(i)] for any values of $y$ such that $Z^2[M(y)]>0$ there exists a continuous derivative, and
\item[(ii)] for any values of $y$ such that $Z^2[M(y)]=0$ we have
\be \label{18}
\lim_{\Delta y\to 0}\frac{M(y+\Delta y)-M(y)}{\Delta y}=+\infty .
\ee
\end{itemize}

Since our function $T=M(y),\ y\geq y_0$ is continuous and increasing (to $+\infty$) there exists unique inverse function that is also continuous and
increasing (to $+\infty$):
\be \label{19}
y=\varphi(T)=\varphi_M(T),\ T\geq T_0[\varphi]=M_\mu(y_0) .
\ee
As a consequence of eqs. (\ref{16}) and (\ref{17}) we have
\be \label{20}
\frac{{\rm d}\varphi(T)}{{\rm d}T}=\frac{Z^2(T)}{\Phi'[\varphi(T)]},\ \Phi'=\Phi'_y[\varphi(T)]>0,\ T\geq T_0[\varphi] ,
\ee
($\varphi'(\gamma)=0$, see (\ref{18})). Eq. (\ref{20}) implies that $\varphi(T)\in C^\infty([T_0[\varphi],+\infty))$ and also that the property (\ref{5})
holds true. Inserting (\ref{19})  into (\ref{14}) one obtains that $X(T)=\varphi_\mu(T),\ T\geq T_0[\varphi]$ is a solution to the integral
equation (\ref{4}).
\end{proof}

\section{Consequences from the TKA formula}

The following Lemma holds true:

\begin{lemma} {\rm
Let $M(y)\in\{ M\}$ be arbitrary, then
\be \label{21}
\int_0^{M(y)} Z^2(t){\rm d}t=F(y)+r(y),
\ee
where
\be \label{22}
F(y)=\frac{y}{2}\ln\left(\frac{y}{2}\right)+E\frac{y}{2}+c_0,\
r(y)=\mcal{O}\left(\frac{\ln(y)}{y}\right),\ E=c-\ln(2\pi) ,
\ee
and the constant within the $\mcal{O}$-symbol is an absolute constant.
}
\end{lemma}

\begin{proof}
We start with the formula (\ref{3}) and $N=1$:
\be \label{23}
\int_0^\infty Z^2(t)e^{-2\delta t}{\rm d}t=\frac{c-\ln(4\pi\delta)}{2\sin(\delta)}+c_0+c_1\delta+\mcal{O}(\delta^2) .
\ee
As long as $|Z(t)|<At^{1/4},\ t\geq t_0$, we have
\bdis
f(t,\delta)=t^{1/2}e^{-\delta t}\leq f\left(\frac{1}{2\delta},\delta\right)=\frac{1}{\sqrt{2 e \delta}},
\edis
\bdis
\int_U^\infty Z^2(t)e^{-2\delta t}{\rm d}t<B\frac{e^{-\delta U}}{\delta^{3/2}},\quad B=\frac{A^2}{\sqrt{2e}},\ U\geq t_0 .
\edis
The value $U=\mu(1/\delta)$ is to be chosen by the following rule:
\bdis
B\delta^{-3/2}e^{-\delta U}\leq \delta^2 \ \Rightarrow \ \mu\left(\frac{1}{\delta}\right)\geq \frac{7}{\delta}\ln\left(\frac{1}{\delta}\right)>
\frac{1}{\delta}\ln\left(\frac{B}{\delta^{7/2}}\right).
\edis
Now, (\ref{23}) implies:
\be \label{24}
\int_0^{\mu(1/\delta)}Z^2(t)e^{-2\delta t}{\rm d}t=\frac{c-\ln(4\pi\delta)}{2\sin(\delta)}+c_0+c_1\delta+\mcal{O}\left(\delta^2\right), \
\mu\left(\frac{1}{\delta}\right)\geq \frac{7}{\delta}\ln\left(\frac{1}{\delta}\right),
\ee
(see the introduction, part (B) - the condition for $\mu(y)$), and for the remainder term we have:
\be \label{25}
-\int_{\mu(1/\delta)}^\infty Z^2(t)e^{-2\delta t}{\rm d}t=\mcal{O}\left(\delta^2\right) .
\ee
Let $\delta\in (0,\delta_0]$ with $\delta_0$ being sufficiently small, then
\bdis
\frac{1}{\sin(\delta)}=\frac{1}{\delta}\left\{ 1+\frac{\delta^2}{6}+\mcal{O}\left(\delta^4\right)\right\},
\edis
and
\bdis
\frac{c-\ln(4\pi\delta)}{2\sin(\delta)}=\frac{1}{2\delta}\ln\left(\frac{1}{\delta}\right)+\frac{D}{2\delta}+
\mcal{O}\left[\delta\ln\left(\frac{1}{\delta}\right)\right],
\edis
and (see (\ref{24}))
\be \label{26}
\int_0^{\mu(1/\delta)}Z^2(t)e^{-2\delta t}{\rm d}t=\frac{1}{2\delta}\ln\left(\frac{1}{\delta}\right)+\frac{D}{2\delta}+c_0+
\mcal{O}\left[\delta\ln\left(\frac{1}{\delta}\right)\right], \ D=c-\ln(4\pi) .
\ee
Putting $\delta=1/y,\ y_0=1/\delta_0$ into eq. (\ref{26}) and using eq. (\ref{14}) we obtain the formulae (\ref{21}) and (\ref{22}),
respectively. Since the constants
in the eqs. (\ref{23}), (\ref{25}) are absolute, the constant entering the $\mcal{O}$-symbol in (\ref{22}) is absolute, too.
\end{proof}

\section{Proof of the theorem}

Putting $T=y/2$ into eq. (\ref{10}) and comparing with the formula (\ref{21}) we obtain
\bdis
\frac{y}{2}<M(y),\ y\to\infty .
\edis
Furthermore, putting into eq. (\ref{10})
\bdis
T=\frac{y}{2}\left(1+\frac{A}{\ln\left(\frac{y}{2}\right)}\right), \ A>1-c ,
\edis
and comparing with eq. (\ref{21}) we have
\bdis
M(y)<\frac{y}{2}\left(1+\frac{A}{\ln\left(\frac{y}{2}\right)}\right),\ y\to\infty.
\edis
Subsequently (see (\ref{19})),
\be \label{27}
0<M(y)-\frac{y}{2}<\frac{A}{2}\frac{y}{\ln\left(\frac{y}{2}\right)}\ \Rightarrow \ 0<2T-\varphi(T)<B\frac{\varphi(T)}{\ln[\varphi(T)]},
\ee
i.e. the following equation holds true
\be \label{28}
1.9 T<\varphi(T)<2T .
\ee
Inserting $y=\varphi(T)\in\{ \varphi\}$ into eq. (\ref{21}) (see (\ref{19})) we obtain the formula (\ref{6}) and the estimate (\ref{7}),
(see (\ref{28})). \\

The relation (\ref{8}) follows from $F'(y)=1/2\ln(y/2)+E+1$ and from the eq. (\ref{6}) written for $\varphi_1$ and $\varphi_2$, respectively, with help
of (\ref{28}). \\

Since $\mu(y)>M(y),\ y\geq y_0$ (see (\ref{14})) and in the case of boundedness of the set $\{\mu(y_0)\}$ the choice of values (see (\ref{9})):
\bdis
T_0=\sup_\mu\{ \mu(y_0)\}\geq \mu(y_0)>M_{\mu}(y_0),\ \forall \mu(y)\in\{\mu\} ,
\edis
is regular, i.e. the interval $[ T_0,+\infty)$ is the common domain of the functions $\varphi_\mu(T),\ \mu(T)\in\{\mu\}$.

\section{Relation between the Jacob's ladders and prime-counting function $\pi(T)$}

Comparing the formulae (\ref{6}) and (\ref{10}) we obtain
\bdis
\omega\left(\frac{\varphi}{2}\right)-\omega(T)=(1-c)T+\mcal{O}\left(T^{1/3+\epsilon}\right), \
\omega(t)=t\ln(t)+(c-\ln(2\pi))T.
\edis
Let us consider the power series expansion in the variable $\varphi/2-T$ of the previous formula. We obtain the following nonlinear equation
\be \label{29}
x(\ln(T)-a)-\sum_{k=2}^\infty \frac{x^k}{k(k-1)}=1-c+\mcal{O}\left(T^{-2/3+\epsilon}\right), \ x=\frac{T-\frac{\varphi}{2}}{T},
\ee
where $a=\ln(2\pi)-1-c$. Because of (see (\ref{27}) and (\ref{28}))
\bdis
x=\mcal{O}\left(\frac{1}{\ln(T)}\right),
\edis
we obtain from (\ref{29}):
\bdis
x=\frac{1-c}{\ln(T)-a}+\mcal{O}\left(\frac{1}{\ln^3(T)}\right)\ \Rightarrow \
T-\frac{\varphi(T)}{2}=\frac{T}{\ln(T)}\left\{ 1+\mcal{O}\left(\frac{1}{\ln(T)}\right)\right\},
\edis
and furthermore, by using the Selberg-Erd\" os theorem, we obtain the formula:
\be \label{30}
\pi(T)\sim \frac{1}{1-c}\left\{ T-\frac{\varphi(T)}{2}\right\},\ T\to\infty ,\  \forall \varphi\in\{\varphi\} .
\ee

\begin{remark} {\rm
As a consequence of the above written we have that the Jacob's ladders are connected (along to the zeroes of the function $\zeta(1/2+it)$) also to the
prime-counting function $\pi(T)$, see (\ref{30}).
}
\end{remark}

More interesting information can be deduced from the nonlinear equation (\ref{29}). Namely, inserting
\be \label{31}
Te^{-a}=\tau,\ e^{-a}\varphi\left( e^a\tau\right)=\psi(\tau)
\ee
into (\ref{29}) we obtain:
\be \label{32}
x\ln(\tau)-\sum_{k=2}^\infty \frac{x^k}{k(k-1)}=1-c+\mcal{O}\left(\tau^{-2/3+\epsilon}\right),\ x=\frac{\tau-\frac{\psi(\tau)}{2}}{\tau} .
\ee
The following statement holds true: if in the equation:
\be \label{33}
x=\frac{A_1}{\ln(\tau)}+\frac{A_3}{\ln^3(\tau)}+\dots +\frac{A_n}{\ln^n(\tau)}+\frac{A_{n+1}}{\ln^{n+1}(\tau)}+\dots \ ,
\ee
the coefficients $A_1,A_3,\dots ,A_n$ are already known, then the coefficient $A_{n+1}$ is determined by (\ref{32}). We obtain:
\begin{eqnarray*}
& & A_1=1-c,\ A_3=\frac{1}{2}(1-c)^2,\ A_4=\frac{1}{6}(1-c)^3, \\
& & A_5=\frac{1}{2}(1-c)^3+\frac{1}{12}(1-c)^4,\ \dots \ .
\end{eqnarray*}

Changing the variables in (\ref{33}) into the initial ones (see (\ref{31})) we obtain the following asymptotic formula
\be \label{34}
\frac{1}{T}\left\{ T-\frac{\varphi(T)}{2}\right\}\sim\frac{A_1}{\ln(T)-a}+\frac{A_3}{(\ln(T)-a)^3}+\dots \sim
\frac{A_1}{\ln(T)}+\frac{B_2}{\ln^2(T)}+\frac{B_3}{\ln^3(T)}+\dots \ , \ T\to\infty ,
\ee
where $B_2=aA_1,\ B_3=a^2A_1+A_3,\ \dots $.

\begin{remark} {\rm
Let us remark that the asymptotic formula (\ref{34}) is an analogue to the following asymptotic formula
\bdis
\frac{1}{T}\int_2^T \frac{{\rm d}t}{\ln(t)}\sim \frac{1}{\ln(T)}+\frac{1}{\ln^2(T)}+\frac{2!}{\ln^3(T)}+\dots \ ,
\edis
for the Gauss logarithmic integral.
}
\end{remark}

\section{Fundamental properties of $Z^2$-transformation}

Let us mention explicitly that the key idea of the proof of the theorem was to introduce the new integral transformation: $Z^2$-transformation. \\

By using the appropriate terminology from optics (see example: Landau\& Lifshitz, Field theory, GIFML, Moscow 1962, page 167):
\begin{itemize}
\item the elements $\mu(y)\in\{\mu\}$ are called the \emph{rays} and the set $\{\mu\}$ itself is called the \emph{beam},

\item the beams crossing each other in a given point are called \emph{homocentric beams}.
\end{itemize}

The fundamental property of the $Z^2$-transformation lies in the following:
\begin{itemize}
\item
if the set $\{ \mu(y_0)\}$ is bounded then the beam $\{\mu\}$
(and, at the same time, also any other homocentric beam) is transformed into the homocentric beam $\{ \varphi\}$
of the Jacob's ladders, with respect to the
point $[+\infty,+\infty]$ (see (\ref{8}),(\ref{9})),
\item the transformed beam $\{\varphi\}$ is very narrow in sense of (\ref{9}), i.e. an $Z^2$-optical system generates an analogue of
a laser beam.
\end{itemize}

Let us consider for example the homocentric (with respect to the point $[y_0,y_0^2]$) sheaf of rays:
\be \label{35}
u(y;\rho,n)=y^2[1+\rho(y-y_0)^n],\ \rho\in[0,1],\ n\in N .
\ee
Following the equation $u_0(y+\Delta;0,n)=u_1(y;1,n)$ we obtain that
\bdis
\Delta=\frac{y(y-y_0)^n}{\sqrt{1+(y-y_0)^2}}\to\infty, \ \mbox{at} \ y\to\infty,
\edis
i.e. (\ref{35}) is a diverging beam. Anyway, the $Z^2$-mapping transforms (\ref{35}) into an analogue of a laser beam $\{\varphi_u\}$.

\section{On intervals that cannot be reached by estimates of Heath-Brown and Ivic}

We will show the accuracy of our formula (2.1) in comparison with known estimates of Heath-Brown and Ivic (see \cite{5}, (7.20) page 178, and
(7.62) page 191)
\be
\int_T^{T+G}Z^2(t){\rm d}t=\mcal{O}\left(G\ln^2(T)\right),\ G=T^{1/3-\epsilon_0},\ \epsilon_0=\frac{1}{108} .
\ee
First of all, one can easily obtain the \emph{tangent law} from both (2.1) and (8.1):
\be
\int_{T}^{T+U}Z^2(t){\rm d}t=U\ln\left(e^{-a}\frac{\varphi(T)}{2}\right)\tan[\alpha(T,U)]+\mcal{O}\left(\frac{1}{T^{1/3+2\epsilon_0}}\right)
\ee
for $0<U<T^{1/3+\epsilon_0}$, where $\alpha=\alpha(U,T)$ is the angle of the chord of the curve $y=\frac{1}{2}\varphi(T)$ crossing the points
$[T,\frac{1}{2}\varphi(T)]$ and $[T+U,\frac{1}{2}\varphi(T+U)]$. Further, from (2.1) and (2.5) we have
\be
\tan(\alpha_0)=\tan[\alpha(T,U_0)]=1+\mcal{O}\left(\frac{1}{\ln(T)}\right),\ U_0=T^{1/3+2\epsilon}.
\ee
And finally, considering the set of all chords of the curve $y=\frac{1}{2}\varphi(T)$ which are parallel to our fundamental chord joining points
$[T,\frac{1}{2}\varphi(T)]$ and $[T+U,\frac{1}{2}\varphi(T+U)]$, we obtain the continuum of formulae:
\begin{eqnarray}
& & \int_{a}^{b}Z^2(t){\rm d}t=(b-a)\ln(T)+\mcal{O}(b-a)+\mcal{O}\left(\frac{1}{T^{1/3+2\epsilon_0}}\right),\nonumber \\
& &
0<b-a<1, (a,b)\subset \left( T,T+T^{1/3+\epsilon_0}\right)   .
\end{eqnarray}

\begin{remark} {\rm
It is quite evident that the interval $\left( 0,1\right)$ cannot be reached in known theories leading to estimates of Heath-Brown and
Ivic.
}
\end{remark}



\begin{thebibliography}{9}
%
%
\bibitem{1}
F. V. Atkinson, `The mean value of the zeta-function on
critical line', {\em
 Quart. J. Math. }10 (1939) 122--128..
%
\bibitem{2}
R. Balasubramanian,
'An improvement on a theorem of Titchmarsh on the mean square of $|\zeta(1/2+it)|$',  {\em Proc. London. Math. Soc.} 3 36 (1978) 540--575.
%
\bibitem{3}
A. Good,
 `Ein $\Omega$-Resultat f\" ur quadratische Mittel der Riemannschen Zetafunktion auf der kritische Linie',
 {\em Invent. Math. } 41 (1977) 233--251.
%
\bibitem{4}
G. H. Hardy and J. E. Littlewood, `Contribution to the theory of the Riemann zeta-function and the theory pf the distribution of Primes',
{\em Acta. Math. } 41 (1918) 119--195.
%
\bibitem{5}
A. Ivic, `The Riemann zeta-function', A Willey-Interscience Publication, New York, 1985.
%
\bibitem{6}
H. Kober, `Eine Mittelwertformeln der Riemannschen Zetafunktion', {\em Composition Math.} 3 (1935) 174--189.
%
\bibitem{7}
E. C. Titchmarsh, `The mean-value of the zeta-function on the critical line', {\em Proc. London. Math. Soc.} (2) 27 (1928) 137--150.
%
\bibitem{8}
E. C. Titchmarsh, `The theory of the Riemann zeta-function', Clarendon Press, Oxford, 1951.
\end{thebibliography}
\end{document}